# The extremogram: A correlogram for extreme events

RICHARD A. DAVIS[1] and THOMAS MIKOSCH[2]

[1]*Department of Statistics, Columbia University, 1255 Amsterdam Avenue, New York, NY 10027, USA. E-mail: rdavis@stat.columbia.edu*
[2]*Laboratory of Actuarial Mathematics, University of Copenhagen, Universitetsparken 5, DK-2100 Copenhagen, Denmark. E-mail: mikosch@math.ku.dk*

We consider a strictly stationary sequence of random vectors whose finite-dimensional distributions are jointly regularly varying with some positive index. This class of processes includes, among others, ARMA processes with regularly varying noise, GARCH processes with normally or Student-distributed noise and stochastic volatility models with regularly varying multiplicative noise. We define an analog of the autocorrelation function, the *extremogram*, which depends only on the extreme values in the sequence. We also propose a natural estimator for the extremogram and study its asymptotic properties under $\alpha$-mixing. We show asymptotic normality, calculate the extremogram for various examples and consider spectral analysis related to the extremogram.

*Keywords:* GARCH; multivariate regular variation; stationary sequence; stochastic volatility process; tail dependence coefficient

## 1. Measures of extremal dependence in a strictly stationary sequence

The motivation for this research comes from the problem of choosing between two popular and commonly used families of models, the generalized autoregressive conditional heteroscedastic (GARCH) process and the heavy-tailed stochastic volatility (SV) process, for modeling a particular financial time series. Both GARCH and SV models possess the *stylized features* exhibited by log-returns of financial assets. Specifically, these time series have heavy-tailed marginal distributions, are dependent but uncorrelated and display stochastic volatility. The latter property is manifested via the often slow decay of the sample autocorrelation function (ACF) of the absolute values and squares of the time series. Since both GARCH and SV models can be chosen to have virtually identical behavior in the tails of the marginal distribution and in the ACF of the squares of the process, it is difficult for a given time series of returns to decide between the two models on the basis of routine time series diagnostic tools.

The problem of finding probabilistically reasonable and statistically estimable measures of extremal dependence in a strictly stationary sequence is, to some extent, an open one. In classical time series analysis, which mostly deals with second order structure of stationary sequences, the ACF is a well-accepted object for describing meaningful information about serial dependence. The ACF is sometimes overvalued as a tool for measuring dependence, especially if one is only interested in extremes. It does, of course, determine the distribution of a stationary Gaussian sequence, but for non-Gaussian and nonlinear time series, the ACF often provides little insight into the dependence structure of the process. This is particularly the case when one considers heavy-tailed nonlinear time series such as the GARCH model. In this case, the estimation of the ACF via the sample ACF is also rather imprecise and even misleading since the asymptotic confidence bands are typically larger than the estimated autocorrelations; see, for example, the results in Basrak *et al.* [1] for bilinear processes, Davis and Mikosch [11], Mikosch and Stărică [26] and Basrak *et al.* [2] for ARCH and GARCH processes and Resnick [32] for teletraffic models.

## 1.1. The extremal index

The asymptotic behavior of the extremes leads to one clear difference between GARCH and SV processes. It was shown in Davis and Mikosch [11], Basrak *et al.* [2], Davis and Mikosch [12] (see also Breidt and Davis [6] for the light-tailed SV case) that GARCH processes exhibit extremal clustering (that is, clustering of extremes), while SV processes lack this form of clustering. Associated with most stationary time series is a parameter $\theta \in (0,1]$, called the *extremal index* (see Leadbetter *et al.* [24]), which is a measure of clustering in the extremes. For example, the extremal index $\theta$ is less than 1 for a GARCH process, which is indicative of extremal clustering, while $\theta = 1$ for SV processes, indicating no clustering. The parameter $\theta$ can also be interpreted as the reciprocal of the expected cluster size in the limiting compound Poisson process of the weakly converging point processes of exceedances of $(X_t)$; see, for example, Leadbetter *et al.* [24] or Embrechts *et al.* [15], Section 8.1.

In this paper, we take a different tack and study the extremal dependence structure of general strictly stationary vector-valued time series $(\mathbf{X}_t)$. Certainly, the cluster distribution of the limiting compound Poisson process contains more useful information about the clustering behavior of extremes than the extremal index. Although explicit formulae for the extremal index and the cluster distribution exist for some specific time series models (including certain ARMA and GARCH models and some Markov processes), these expressions are, in general, very complicated to compute and even difficult to simulate. They are also rather difficult objects to estimate and do not always yield satisfactory results, even for moderate sample sizes.

## 1.2. Regularly varying time series

In this paper, we focus on strictly stationary sequences whose finite-dimensional distributions have power law tails in some generalized sense. In particular, we will assume that the finite-dimensional distributions of the $d$-dimensional process $(\mathbf{X}_t)$ have regularly varying distributions with a positive tail index $\alpha$. This means that for any $h \geq 1$, the lagged vector $\mathbf{Y}_h = \mathrm{vec}(\mathbf{X}_1, \ldots, \mathbf{X}_h)$ satisfies the relation

$$\frac{P(x^{-1}\mathbf{Y}_h \in \cdot)}{P(|\mathbf{Y}_h| > x)} \stackrel{v}{\to} \mu_h(\cdot) \tag{1.1}$$

for some non-null Radon measure $\mu_h$ on $\overline{\mathbb{R}}_0^{hd} = \overline{\mathbb{R}}^{hd} \setminus \{\mathbf{0}\}$, $\overline{\mathbb{R}} = \mathbb{R} \cup \{\pm\infty\}$, with the property that $\mu_h(tC) = t^{-\alpha}\mu_h(C)$, $t > 0$, for any Borel set $C \subset \overline{\mathbb{R}}_0^{hd}$. Here, $\stackrel{v}{\to}$ denotes vague convergence; see Kallenberg [22], Daley and Vere-Jones [9] and Resnick [30, 31, 33] for this notion and Resnick [30, 31, 33] and Hult and Lindskog [21] for the notion of multivariate regular variation. We call such a sequence $(\mathbf{X}_t)$ *regularly varying with index* $\alpha > 0$.

Various time series models of interest are regularly varying. These include infinite variance stable processes, ARMA processes with i.i.d. regularly varying noise, GARCH processes with i.i.d. noise with infinite support (including normally and Student-distributed noise) and stochastic volatility models with i.i.d. regularly varying noise. In Section 2, we will be more precise about the regular variation of the aforementioned sequences. It follows from general multivariate extreme value theory (e.g., Resnick [31]) that any strictly stationary time series whose finite-dimensional distributions are in the maximum domain of attraction of a multivariate extreme value distribution can be transformed to a regularly varying strictly stationary time series. This can be simply achieved by a monotone transformation of the marginal distribution. Hence, the results of this paper apply in a more general framework than that of regularly varying sequences.

For our purposes, it will be convenient to use a sequential definition of a regularly varying sequence $(\mathbf{X}_t)$ which is equivalent to the definition above. Consider a sequence $a_n \uparrow \infty$ such that $P(|\mathbf{X}| > a_n) \sim n^{-1}$. Then, (1.1) holds if and only if there exist constants $c_h > 0$ such that

$$nP(a_n^{-1}\mathbf{Y}_h \in \cdot) \stackrel{v}{\to} c_h\mu_h(\cdot) = \nu_h(\cdot), \tag{1.2}$$

where $\mu_h$ is defined in (1.1). Alternatively, for each $h \geq 1$, one can replace $(a_n)$ in (1.2) by a sequence $(a_n^{(h)})$ such that $P(|\mathbf{Y}_h| > a_n^{(h)}) \sim n^{-1}$ and then $c_h = 1$ in (1.2). However, for each $h \geq 1$, $a_n/a_n^{(h)} \to d_h$ as $n \to \infty$ for some positive constants $d_h$.

## 1.3. The upper tail dependence coefficient

As a starting point for the definition of a measure of extremal dependence in a strictly stationary sequence, we consider the (upper) *tail dependence coefficient*. It is defined for a two-dimensional vector $(X,Y)$ with $X\stackrel{d}{=}Y$ as the limit (provided it exists)

$$\lambda(X,Y)=\lim_{x\to\infty}P(X>x|Y>x).$$

Of course, $\lambda\in[0,1]$, and $\lambda=0$ when $X$ and $Y$ are independent or asymptotically independent. The larger the $\lambda$, the larger the extremal dependence in the vector $(X,Y)$. We refer, for example, to the discussions in Ledford and Tawn [25] and Beirlant *et al.* [3] on the tail dependence coefficient.

The tail dependence coefficient can also be applied to the pairs $(X_0,X_h)$ of a one-dimensional strictly stationary time series. The collection of values $\lambda(X_0,X_h)$ contains useful information about the serial extremal dependence in the sequence $(X_t)$. If one considers a real-valued, regularly varying sequence $(X_t)$ with index $\alpha>0$, the definition of regular variation immediately ensures the existence of the quantities

$$\begin{aligned}\lambda(X_0,X_h)&=\lim_{x\to\infty}\frac{P(x^{-1}(X_0,\dots,X_h)\in(1,\infty)\times(0,\infty)^{h-1}\times(1,\infty))}{P(x^{-1}(X_0,\dots,X_h)\in(1,\infty)\times(0,\infty)^{h})}\\&=\frac{\mu_{h+1}((1,\infty)\times(0,\infty)^{h-1}\times(1,\infty))}{\mu_{h+1}((1,\infty)\times(0,\infty)^{h})}.\end{aligned}$$

## 1.4. The extremogram

Now, let $(\mathbf{X}_t)$ be a strictly stationary, regularly varying sequence of $\mathbb{R}^d$-valued random vectors. Consider two Borel sets $A,B$ in $\overline{\mathbb{R}}^d$ such that $C=A\times\overline{\mathbb{R}}^{d(h-1)}\times B$ is bounded away from zero and $\nu_{h+1}(\partial C)=0$. According to (1.2), the following limit exists:

$$nP(a_n^{-1}\mathbf{X}_0\in A,a_n^{-1}\mathbf{X}_h\in B)\to\nu_{h+1}(A\times\overline{\mathbb{R}}_0^{d(h-1)}\times B)=\gamma_{AB}(h).$$

Note that if both $A$ and $B$ are bounded away from zero, then

$$\begin{aligned}n\operatorname{cov}(I_{\{a_n^{-1}\mathbf{X}_0\in A\}},I_{\{a_n^{-1}\mathbf{X}_h\in B\}})&=n[P(a_n^{-1}\mathbf{X}_0\in A,a_n^{-1}\mathbf{X}_h\in B)-P(a_n^{-1}\mathbf{X}\in A)P(a_n^{-1}\mathbf{X}\in B)]\\&\sim nP(a_n^{-1}\mathbf{X}_0\in A,a_n^{-1}\mathbf{X}_h\in B)\sim\gamma_{AB}(h).\end{aligned}$$

Also, note that the strictly stationary bivariate time series $(I_{\{a_n^{-1}\mathbf{X}_t\in A\}},I_{\{a_n^{-1}\mathbf{X}_t\in B\}})$ has limiting covariance matrix function given by

$$\mathbf{\Gamma}(h)=\begin{bmatrix}\gamma_{AA}(h)&\gamma_{AB}(h)\\\gamma_{BA}(h)&\gamma_{BB}(h)\end{bmatrix},$$

which has all non-negative components. Since $\mathbf{\Gamma}(h)$ is the limit of a sequence of covariance matrix functions, it must also be a covariance matrix function and hence a non-negative definite matrix-valued function; see Brockwell and Davis [8]. In particular, both $(\gamma_{AA}(h))$ and $(\gamma_{BB}(h))$ are non-negative definite functions and $(\gamma_{AB}(h))$ is a cross-covariance function and need not be symmetric in $A$ and $B$.

Alternatively, for $A$ and $A \times B$ bounded away from zero and with $\nu_1(A) > 0$, one may consider

$$P(a_n^{-1}\mathbf{X}_h \in B | a_n^{-1}\mathbf{X}_0 \in A) = \frac{P(a_n^{-1}\mathbf{X}_0 \in A, a_n^{-1}\mathbf{X}_h \in B)}{P(a_n^{-1}\mathbf{X} \in A)} \to \frac{\gamma_{AB}(h)}{\nu_1(A)} = \rho_{AB}(h).$$

Then, since

$$\frac{\operatorname{cov}(I_{\{a_n^{-1}\mathbf{X}_0 \in A\}}, I_{\{a_n^{-1}\mathbf{X}_h \in A\}})}{\operatorname{var}(I_{\{a_n^{-1}\mathbf{X} \in A\}})} \sim \rho_{AA}(h),$$

$(\rho_{AA}(h))$ is the correlation function of a stationary process. With the exception of $A = B$, $(\rho_{AB}(h))$ and the correlation function (with $\nu_1(A)\nu_1(B) > 0$)

$$\operatorname{corr}(I_{\{a_n^{-1}\mathbf{X}_0 \in A\}}, I_{\{a_n^{-1}\mathbf{X}_h \in B\}}) \sim \frac{\gamma_{AB}(h)}{\sqrt{\nu_1(A)\nu_1(B)}}, \qquad h \in \mathbb{Z},$$

are, in general, different functions. However, for fixed $A$, all of these quantities are proportional to each other. In what follows, we refer to any one of these limiting quantities, considered as a function of $h$, as the *extremogram* of the sequence $(X_t)$. Since $A, B$ are arbitrary, there exist infinitely many extremograms. The sequence of the tail dependence coefficients of a regularly varying one-dimensional strictly stationary sequence $(X_t)$ is a special case of the extremogram. Indeed,

$$\lambda(X_0, X_h) = \rho_{(1,\infty),(1,\infty)}(h).$$

As mentioned above, it can be interpreted as a particular ACF.

Since $\gamma_{AA}$ can be interpreted as an autocovariance function, one can translate various notions from classical time series analysis to the extremogram. For example, one can introduce the analog of the spectral distribution corresponding to $\gamma_{AA}$. In particular, if $\gamma_{AA}$ is summable, then there exists a spectral density and one may speak of short range dependence in the time series context. Alternatively, if $\gamma_{AA}$ is not summable, then one can talk of long range dependence.

The bivariate extremal dependence measure $\gamma_{AB}$ introduced above can be extended in such a way that any finite number of events $A_1, \ldots, A_h$ is involved. Provided the set $C = A_1 \times \cdots \times A_h$ is bounded away from zero in $\overline{\mathbb{R}}_0^{dh}$ and $\nu_h(\partial C) = 0$, one can define the limiting dependence measure

$$nP(a_n^{-1}\mathbf{X}_1 \in A_1, \ldots, a_n^{-1}\mathbf{X}_h \in A_h) \to \nu_h(A_1 \times \cdots \times A_h).$$

Such quantities can be of interest, for example, when considering the limits of conditional probabilities of the form

$$P(a_n^{-1}\mathbf{X}_2 \in A_2, \ldots, a_n^{-1}\mathbf{X}_h \in A_h | a_n^{-1}\mathbf{X}_1 \in A_1), \tag{1.3}$$

where $A_1$ is bounded away from zero. Probabilities of this form and their limits appear as the extremal index and the cluster probability distribution of strictly stationary sequences; see also Fasen *et al.* [16] who consider a generalization of the tail dependence coefficient. In this paper, we focus on the two-dimensional case, that is, the extremogram, but, in a sense, the extremogram also covers this more general case. Indeed, if we define the strictly stationary process $\mathbf{Y}_h = \operatorname{vec}(\mathbf{X}_1, \dots, \mathbf{X}_h)$, then (1.3) can be written in the form

$$P(a_n^{-1}\mathbf{Y}_h \in \overline{\mathbb{R}}^d \times A_2 \times \cdots \times A_h | \mathbf{Y}_1 \in \overline{\mathbb{R}}^{d(h-1)} \times A_1)$$

whose limit is an extremogram.

The paper is organized as follows. In Section 2, we consider some well-known time series models, including the GARCH and SV models, and discuss conditions under which they constitute a regularly varying sequence. The extremograms are also computed for these models. In Section 3, we study estimators of the extremogram. Assuming that the sequence $(\mathbf{X}_t)$ meets certain dependence conditions such as $\alpha$-mixing with a suitable rate, we show that these estimators are asymptotically unbiased and satisfy a central limit theorem. In Section 4, we apply the asymptotic results to GARCH and SV models. The Fourier transform of the extremogram can be viewed as the analog of the spectral density of a correlogram. The periodogram is similarly defined as the Fourier transform of the estimated extremogram. In Section 5, we show that the periodogram is asymptotically unbiased for the spectral density. A lag window estimate of the spectral density is also formulated and shown to be asymptotically unbiased and consistent. The proof of the main theorem in Section 3 is provided in Section 6.

## 2. Examples of extremograms

### 2.1. Preliminaries on regular variation

We will often make use of a multivariate version of a result of Breiman [7] which can be found in Basrak *et al.* [2]. Assume that the $d$-dimensional vector $\mathbf{X}$ is regularly varying with index $\alpha$ and limiting measure $\mu$, that is, $\mathbf{Y}_h$ and $\mu_h$ in (1.1) are replaced by $\mathbf{X}$ and $\mu$, respectively. Let $\mathbf{A}$ be a random $k \times d$ matrix that is independent of $\mathbf{X}$ with $E\|\mathbf{A}\|^{\alpha+\epsilon} < \infty$ for some $\epsilon > 0$. Then,

$$\frac{P(x^{-1}\mathbf{A}\mathbf{X} \in \cdot)}{P(|\mathbf{X}| > x)} \stackrel{v}{\to} E\mu(\{\mathbf{x} \in \overline{\mathbb{R}}_0^d \colon \mathbf{A}\mathbf{x} \in \cdot\}), \qquad x \to \infty, \tag{2.1}$$

where $\stackrel{v}{\to}$ refers to vague convergence in $\overline{\mathbb{R}}_0^k$.

### 2.2. The stochastic volatility model

We consider a stochastic volatility model

$$X_t = \sigma_t Z_t,$$

where the volatility sequence $(\sigma_t)$ constitutes a strictly stationary sequence of non-negative random variables, independent of the i.i.d. sequence $(Z_t)$. We further assume that $Z$ is regularly varying with index $\alpha > 0$, that is, the limits

$$p = \lim_{x\to\infty} \frac{P(Z > x)}{P(|Z| > x)} \quad \text{and} \quad q = \lim_{x\to\infty} \frac{P(Z \le -x)}{P(|Z| > x)}$$

exist for some $p, q \ge 0$ with $p + q = 1$ and $P(|Z| > x) = x^{-\alpha} L(x)$ for some slowly varying function $L$. If we also assume that $E(\sigma^{\alpha+\epsilon}) < \infty$ for some $\epsilon > 0$, then $(X_t)$ is regularly varying with index $\alpha$. This follows from Breiman's result (2.1); see Davis and Mikosch [12]. Hence, the finite-dimensional distributions of $(X_t)$ are regularly varying with index $\alpha$ and the limiting measures $\nu_h$ in (1.2) are concentrated at the axes, as in the case of an i.i.d. sequence. Equivalently, the corresponding spectral measures are concentrated at the intersection of the unit sphere with the axes. To be precise (see [12]),

$$\nu_{h+1}(dx_0, \ldots, dx_h) = \sum_{i=0}^{h} \lambda_\alpha(dx_i) \prod_{i\ne j} \varepsilon_0(dx_j),$$

where $\lambda_\alpha(x, \infty] = px^{-\alpha}$ and $\lambda_\alpha[-\infty, -x] = qx^{-\alpha}$, $x > 0$, and $\varepsilon_y$ denotes Dirac measure at $y$. In particular, if $h \ge 1$ and $C = A \times \overline{\mathbb{R}}^{h-1} \times B$ is bounded away from zero with $\nu_{h+1}(\partial C) = 0$, then

$$\gamma_{AB}(h) = \begin{cases} \lambda_\alpha(A), & \text{if } 0 \in B,\ A \text{ is bounded away from zero,} \\ \lambda_\alpha(B), & \text{if } 0 \in A,\ B \text{ is bounded away from zero,} \\ 0, & \text{if } A \text{ and } B \text{ are both bounded away from zero.} \end{cases}$$

This means that $\gamma_{AB}(h)$ is zero unless either $A$ or $B$ contains zero. Moreover, $\rho_{AB}(h) = 1$ if $0 \in B$ and $\rho_{AB}(h) = 0$ otherwise. In particular, $\gamma_{AB}(h)$ does not depend on $h$ for $h \ge 1$. If both $A$ and $B$ contain 0, then the set $A \times \overline{\mathbb{R}}_0^{h-1} \times B$ is not bounded away from zero.

## 2.3. GARCH process

The regular variation of a GARCH$(p, q)$ process was shown in Basrak *et al.* [2] under general conditions. Here, we focus on a GARCH$(1, 1)$ process because the calculations can be made explicit; this is not possible for a general GARCH process. A GARCH$(1, 1)$ process is given by the equations

$$X_t = \sigma_t Z_t, \qquad t \in \mathbb{Z},$$

where $(Z_t)$ is an i.i.d. sequence with $EZ = 0$ and $\mathrm{var}(Z) = 1$, and

$$\sigma_t^2 = \alpha_0 + \alpha_1 X_{t-1}^2 + \beta_1 \sigma_{t-1}^2 = \alpha_0 + \sigma_{t-1}^2 C_{t-1}, \qquad C_t = \alpha_1 Z_t^2 + \beta_1. \tag{2.2}$$

The parameters $\alpha_0, \alpha_1, \beta_1 > 0$ are chosen such that $(X_t)$ is strictly stationary and the unique positive solution to the equation

$$EC^{\alpha/2} = 1 \tag{2.3}$$

exists. Then, under regularity conditions such as the existence of a positive density of $Z$ on $\mathbb{R}$, the sequences $(\sigma_t)$ and $(X_t)$ are regularly varying with index $\alpha$. This follows from theory developed by Kesten [23]. Equation (2.3) has a positive solution if $Z$ is standard normal or Student distributed. We refer to Mikosch and Stărică [26], Theorem 2.6, for details in the GARCH$(1,1)$ case.

We now calculate the extremogram $\gamma_{AB}$ for the sets $A = (a, \infty)$, $B = (b, \infty)$ for positive $a, b$. For more general sets, the calculations become less tractable. We will make repeated use of the following auxiliary result whose proof can be found in Mikosch and Stărică [26] and Basrak *et al.* [2].

**Lemma 2.1.** *Assume that the strictly stationary* GARCH$(1,1)$ *process* $(X_t)$ *satisfies the conditions of Theorem* 2.6 *in Mikosch and Stărică [26]. Then,* $(X_t)$ *is regularly varying with index* $\alpha > 0$ *given as the solution to (2.3) and the following relations hold for any* $h \geq 2$:

$$(\sigma_1^2, \ldots, \sigma_h^2) = \sigma_0^2 (C_0, C_1 C_0, \ldots, C_{h-1} \cdots C_0) + \mathbf{R}_1, \tag{2.4}$$

$$(X_1^2, \ldots, X_h^2) = \sigma_1^2 (Z_1^2, Z_2^2 C_1, \ldots, Z_h^2 C_{h-1} \cdots C_1) + \mathbf{R}_2, \tag{2.5}$$

*where, for any* $\epsilon > 0$, $nP(n^{-2/\alpha} |\mathbf{R}_i| > \epsilon) \to 0, i = 1, 2.$

It follows from Lemma 2.1 and Breiman's result (2.1) that

$$\begin{aligned} \frac{P(x^{-1} X_0^2 \in A, x^{-1} X_h^2 \in B)}{P(x^{-1} X^2 \in A)} &\sim \frac{P(x^{-1} \sigma_0^2 Z_0^2 / a > 1, x^{-1} \sigma_0^2 C_0 \cdots C_{h-1} Z_h^2 / b > 1)}{P(X^2 / a > x)} \\ &\sim \frac{E(\min(Z_0^2 / a, C_0 \cdots C_{h-1} Z_h^2 / b))^{\alpha/2}}{E(Z^2 / a)^{\alpha/2}} = \rho_{AB}(h). \end{aligned}$$

It is, in general, not possible to obtain more explicit expressions for $\rho_{AB}$. In the ARCH$(1)$ case, that is, when $\beta_1 = 0$, we can use (2.3) to obtain

$$\rho_{AB}(h) = \frac{E(\min(C_0, C_0 \cdots C_h a / b))^{\alpha/2}}{EC^{\alpha/2}} = E(\min(1, C_0 \cdots C_{h-1} a / b))^{\alpha/2}.$$

The right-hand side decays to zero at an exponential rate. This can also be seen from the following calculations in the GARCH$(1,1)$ case. There exists some constant $c > 0$ such that

$$\begin{aligned} \rho_{AB}(h) &\leq cE(\min(\alpha_1 Z_0^2, C_0 \cdots C_{h-1}(\alpha_1 Z_h^2)))^{\alpha/2} \\ &\leq cE(\min(C_0, C_0 \cdots C_h))^{\alpha/2} \end{aligned}$$

$$
\begin{aligned}
&= cE(\min(1, C_0 \cdots C_{h-1}))^{\alpha/2} \\
&= cP(C_0 \cdots C_{h-1} \geq 1) + E(C_0 \cdots C_{h-1})^{\alpha/2} I_{\{C_0 \cdots C_{h-1} < 1\}}.
\end{aligned}
$$

Choose $\kappa \in (0, \alpha/2)$. Since the function $r \to EC^r$ is convex and (2.3) holds, $EC^\kappa < 1$. Then, by Markov's inequality, $P(C_0 \cdots C_{h-1} \geq 1) \leq (EC^\kappa)^h$ and

$$
E((C_0 \cdots C_{h-1})^{\alpha/2} I_{\{C_0 \cdots C_{h-1} < 1\}}) \leq E(C_0 \cdots C_{h-1})^\kappa = (EC^\kappa)^h.
$$

Hence, $\rho_{AB}(h) \leq c(EC^\kappa)^h$, which implies that $\rho_{AB}(h)$ decays to zero exponentially fast. In [13], we give some further examples of extremograms for a GARCH$(1,1)$ process.

## 2.4. Symmetric $\alpha$-stable processes

Let $(X_t)$ be a strictly stationary symmetric $\alpha$-stable (s$\alpha$s) sequence with integral representation

$$
X_t = \int_E f_t \, \mathrm{d}M, \qquad t \in \mathbb{Z}, \tag{2.6}
$$

where $(f_t)$ is a sequence of deterministic functions such that $f_t \in L^\alpha(E, \mathcal{E}, m)$ for some $\alpha \in (0, 2)$, $\mathcal{E}$ is a $\sigma$-field on $E$ and $m$ is a measure on $\mathcal{E}$. The measure $m$ is the control measure of the s$\alpha$s random measure $M$ on $E$. For the definition of $\alpha$-stable integrals of type (2.6), we refer to Samorodnitsky and Taqqu [37]. Conditions for stationarity of the sequence $(X_t)$ were given by Rosiński [34]. By the definition of an s$\alpha$s integral and stationarity of $(X_t)$, for some constant $C_\alpha > 0$, the tail of the marginal distribution satisfies

$$
P(X_t > x) \sim C_\alpha x^{-\alpha} \int_E |f_t|^\alpha \, \mathrm{d}m = C_\alpha x^{-\alpha} \int_E |f_0|^\alpha \, \mathrm{d}m.
$$

The next result follows from Samorodnitsky [36]; also see Theorem 3.5.6 in Samorodnitsky and Taqqu [37] or Theorem 8.8.18 in Embrechts *et al.* [15]. We have, for $A = (a, \infty)$, $B = (b, \infty)$, $a, b > 0$,

$$
\begin{aligned}
\frac{P(x^{-1} X_h > b, x^{-1} X_0 > a)}{P(x^{-1} X > a)} &\sim \frac{P(x^{-1} \min(X_0/a, X_h/b) > 1)}{P(x^{-1} X > a)} \\
&\sim \frac{\int_E [(\min(f_0^+, f_h^+(a/b)))^\alpha + (\min(f_0^-, f_h^-(a/b)))^\alpha] \, \mathrm{d}m}{\int_E |f_0|^\alpha \, \mathrm{d}m} \\
&= \rho_{AB}(h).
\end{aligned}
$$

If we choose $E = \mathbb{R}$, $m$ to be Lebesgue measure on $\mathbb{R}$ and

$$
f_t(x) = \mathrm{e}^{-\lambda(t-x)} I_{(-\infty, t]}(x),
$$

then the corresponding process $(X_t)_{t\in\mathbb{Z}}$ is the discrete version of an s$\alpha$s Ornstein–Uhlenbeck process and

$$\rho_{AB}(h) = \frac{\int_{-\infty}^{0} \mathrm{e}^{\lambda\alpha x} \min(1, \mathrm{e}^{-\lambda\alpha h}(a/b)^{\alpha})\,\mathrm{d}x}{\int_{-\infty}^{0} \mathrm{e}^{\lambda\alpha x}\,\mathrm{d}x} = \min(1, \mathrm{e}^{-\lambda\alpha h}(a/b)^{\alpha}), \qquad h \ge 0.$$

For $\alpha = 2$ and $a/b = 1$, this autocorrelation function coincides with the autocorrelation function of a Gaussian AR(1) process.

If we assume that

$$f_t(x) = f(t-x), \qquad x \in \mathbb{R}, t \in \mathbb{Z},$$

and $f$ is constant on the intervals $(n-1, n]$ for all $n \in \mathbb{Z}$, then $(X_t)$ is a linear process with i.i.d. s$\alpha$s noise. In this case,

$$\rho_{AB}(h) = \sum_{n=-\infty}^{\infty} ((\min([f(n)]^+, [f(h+n)]^+(a/b)))^{\alpha} + (\min([f(n)]^-, [f(h+n)]^-(a/b)))^{\alpha}) \Big/ \sum_{n=-\infty}^{\infty} |f(n)|^{\alpha}.$$

## 2.5. ARMA process

The extremogram for an ARMA process generated by heavy-tailed noise can be derived directly from the previous example. Suppose that $(X_t)$ satisfies the ARMA$(p,q)$ recursions

$$X_t = \phi_1 X_{t-1} + \cdots + \phi_p X_{t-P} + Z_t + \theta_1 Z_{t-1} + \cdots + \theta_q Z_{t-q},$$

where the autoregressive polynomial $\phi(z) = 1 - \phi_1 z - \cdots - \phi_p z^p$ has no zeros inside or on the unit circle and $(Z_t)$ is an i.i.d. sequence of symmetric and regularly varying random variables. Then, $(X_t)$ has the causal representation

$$X_t = \sum_{j=0}^{\infty} \psi_j Z_{t-j},$$

where the coefficients $\psi_j$ are found from the relation $\phi(z)\sum_{j=0}^{\infty}\psi_j z^j = (1 + \theta_1 z + \cdots + \theta_q z^q)$ (see Brockwell and Davis [8]). From the previous example with the same sets $A$ and $B$, the extremogram is given by

$$\rho_{AB}(h) = \frac{\sum_{j=0}^{\infty}((\min(\psi_j^+, \psi_{j+h}^+(a/b)))^{\alpha} + (\min(\psi_j^-, \psi_{j+h}^-(a/b)))^{\alpha})}{\sum_{n=0}^{\infty}|\psi_j|^{\alpha}}, \qquad h \ge 0.$$

In particular, if $(X_t)$ is an AR(1) process with $\phi_1 \in (0,1)$, then $\psi_j = \phi_1^j$ and

$$\rho_{AB}(h) = \min(1, \phi_1^{\alpha h}(a/b)^{\alpha}).$$

Regardless of the values of $a$ and $b$, the extremogram eventually decays at a geometric rate. It is worth noting that for the case $a > b$, the extremogram may be equal to one for several lags before beginning its exponential descent. If we assume that $a = b = 1$ and $\phi \in (-1, 0)$, then we get

$$\rho_{AA}(2h+1) = 0 \quad \text{and} \quad \rho_{AA}(2h) = |\phi|^{\alpha 2h}.$$

This means that an AR(1) process with a negative coefficient has as an alternating extremogram $\rho_{AA}$ that is zero for all odd lags and decays geometrically for even lags. In this case, the extremogram coincides with the ACF of an AR(2) process with lag-1 coefficient equal to 0 and lag-2 coefficient equal to $|\phi|^{\alpha 2}$. Based on the empirical estimate of the extremogram, AR-type behavior with a negative parameter $\phi$ can be observed for foreign exchange rate high frequency data. See Figure 1 for an illustration.

## 3. Consistency and a central limit theory for the empirical extremogram

The aim of this section is to derive relevant asymptotics for the empirical extremogram. In Sections 3.1 and 3.2, we establish key large-sample properties for the empirical estimator of $\mu(C)$. Based on these results, the asymptotic normality for the empirical extremogram is established in Section 3.3.

Throughout this section, it is assumed that $(\mathbf{X}_t)$ is a strictly stationary, regularly varying sequence with index $\alpha > 0$. The vector $\mathbf{X} = \mathbf{X}_0$ assumes values in $\mathbb{R}^d$ and has limiting measure $\mu$. This means that we replace $\mathbf{Y}_h$ with $\mathbf{X}$ and $\mu_h$ with $\mu$ in the definition of (1.1).

The empirical extremogram, defined in Section 3.3, can be viewed as a ratio of estimates of $\mu(A)$ and $\mu(B)$ for two suitably chosen sets $A$ and $B$. We first consider estimates of $\mu(C)$, where $C$ is a generic subset of $\overline{\mathbb{R}}_0^d$, bounded away from zero and with $\mu(\partial C) = 0$. Then, in particular,

$$\text{there exists} \qquad \epsilon > 0 \qquad \text{such that} \qquad C \subset \{\mathbf{x} \in \overline{\mathbb{R}}^d \colon\ |\mathbf{x}| > \epsilon\}. \tag{3.1}$$

A natural estimator of $\mu(C)$ is given by

$$\widehat{P}_m(C) = \frac{m_n}{n} \sum_{t=1}^{n} I_{\{\mathbf{X}_t/a_m \in C\}},$$

where $(a_n)$ is chosen such that $P(|\mathbf{X}| > a_n) \sim n^{-1}$, $m = m_n \to \infty$ and $m_n/n = o(1)$. These conditions on $(m_n)$ ensure consistency of $\widehat{P}_m(C)$; see Theorem 3.1. The estimator $\widehat{P}_m(C)$ is closely related to the tail empirical process. We refer to the recent monographs de Haan and Ferreira [19], Resnick [33] and the references cited therein.

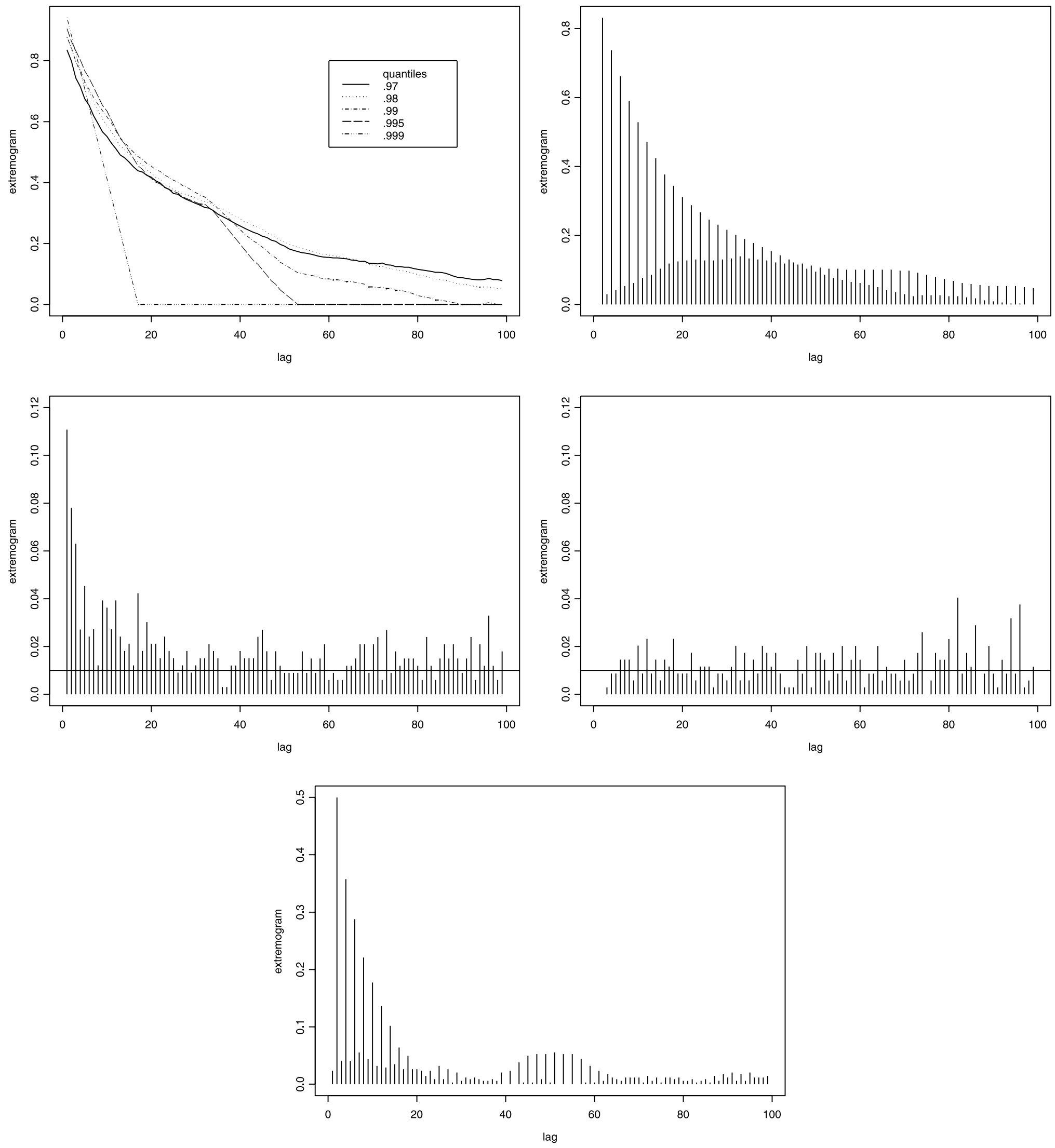

**Figure 1.** The empirical extremogram with $A = B = (1, \infty) \cup (\infty, -1)$ (upper-left) and $A = B = (1, \infty)$ (upper-right) of a sample of five-minute return data of the foreign exchange rate USD-DEM. The middle consists of the extremogram of the residuals from a fitted AR(18) model and from the residuals of a GARCH fitted to the residuals $A = B = (1, \infty)$. The bottom shows the extremogram of a simulated AR-GARCH model. See Section 3.4 for more details.

We will work under the following mixing/dependence conditions on the sequence $(\mathbf{X}_t)$:

(M) The sequence $(\mathbf{X}_t)$ is $\alpha$-mixing with rate function $(\alpha_t)$. Moreover, there exist $m_n$, $r_n \to \infty$ with $m_n/n \to 0$ and $r_n/m_n \to 0$ such that

$$\lim_{n\to\infty} m_n \sum_{h=r_n}^{\infty} \alpha_h = 0 \tag{3.2}$$

and, for all $\epsilon > 0$,

$$\lim_{k\to\infty} \limsup_{n\to\infty} m_n \sum_{h=k}^{r_n} P(|\mathbf{X}_h| > \epsilon a_m, |\mathbf{X}_0| > \epsilon a_m) = 0. \tag{3.3}$$

Condition (3.3) is similar in spirit to condition (2.8) used in Davis and Hsing [10] for establishing convergence of a sequence of point processes to a limiting cluster point process. It is much weaker than the anti-clustering condition $D'(\epsilon a_n)$ of Leadbetter, which is well known in the extreme value literature; see Leadbetter *et al.* [24] or Embrechts *et al.* [15]. Condition (3.3) is equivalent to

$$\lim_{k\to\infty} \limsup_{n\to\infty} \sum_{h=k}^{r_n} P(|\mathbf{X}_h| > \epsilon a_m \,|\, |\mathbf{X}_0| > \epsilon a_m) = 0. \tag{3.4}$$

There are various time series models that are $\alpha$-mixing with geometric rate and for which (3.2) and (3.3) are easily verified. These include GARCH, stochastic volatility and ARMA models under suitable conditions on the noise; see the discussion in Sections 4.1 and 4.2 for GARCH and SV models.

## 3.1. Asymptotic mean and variance

In this section, we calculate the asymptotic mean and variance of $\widehat{P}_m(C)$ under condition (M).

**Theorem 3.1.** *Assume that $(\mathbf{X}_t)$ is a regularly varying, strictly stationary $\mathbb{R}^d$-valued sequence with index $\alpha > 0$ in the sense of (1.1). Moreover, let $C$ and $C \times \overline{\mathbb{R}}_0^{d(h-1)} \times C \subset \overline{\mathbb{R}}_0^{(h+1)d}$ be continuity sets with respect to $\mu$ and $\mu_{h+1}$ for $h \ge 1$, respectively, and let $C$ be bounded away from zero. If condition (M) holds, then*

$$E\widehat{P}_m(C) \to \mu(C), \tag{3.5}$$

$$\operatorname{var}(\widehat{P}_m(C)) \sim \frac{m_n}{n}\left[\mu(C) + 2\sum_{h=1}^{\infty} \tau_h(C)\right], \tag{3.6}$$

*where $\tau_h(C) = \mu_{h+1}(C \times \overline{\mathbb{R}}_0^{d(h-1)} \times C)$. If $\mu(C) = 0$, then (3.6) is interpreted as $\operatorname{var}(\widehat{P}_m(C)) = o(m_n/n)$. In particular, we have $\widehat{P}_m(C) \xrightarrow{P} \mu(C)$.*

**Proof.** In what follows, it will be convenient to write

$$P_m(C) = mP(\mathbf{X}/a_m \in C) = mp_0 \quad \text{and} \quad p_{st} = P(\mathbf{X}_s/a_m \in C, \mathbf{X}_t/a_m \in C).$$

Regular variation of $\mathbf{X}$ and strict stationarity of $(\mathbf{X}_t)$ imply that

$$E\widehat{P}_m(C) = P_m(C) \to \mu(C) \qquad \text{as } n \to \infty.$$

This proves (3.5).

Turning to (3.6), we first note that

$$\begin{aligned} \operatorname{var}(\widehat{P}_m(C)) &= \left(\frac{m_n}{n}\right)^2 n \operatorname{var}(I_{\{\mathbf{X}/a_m \in C\}}) + 2\left(\frac{m_n}{n}\right)^2 \sum_{h=1}^{n-1} (n-h) \operatorname{cov}(I_{\{\mathbf{X}_0/a_m \in C\}}, I_{\{\mathbf{X}_h/a_m \in C\}}) \\ &= I_1 + I_2. \end{aligned}$$

By regular variation of $\mathbf{X}$,

$$I_1 = \frac{m_n}{n}[P_m(C)(1-p_0)] \sim \frac{m_n}{n}\mu(C). \tag{3.7}$$

We have, for $k \geq 1$ fixed,

$$\begin{aligned} \frac{n}{2m_n} I_2 &= m_n \left( \sum_{h=1}^{k} + \sum_{h=k+1}^{r_n} + \sum_{h=r_n+1}^{n-1} \right) (1-h/n)[p_{0h} - p_0^2] \\ &= I_{21} + I_{22} + I_{23}. \end{aligned}$$

The regular variation of $(\mathbf{X}_t)$ implies that $I_{21} \to \sum_{i=1}^{k} \tau_h(C)$, so it suffices to show that

$$\lim_{k\to\infty} \limsup_{n\to\infty} (|I_{22}| + |I_{23}|) = 0. \tag{3.8}$$

Since $C$ is bounded away from zero, (3.1) holds. Then, since $r_n = \mathrm{o}(m_n)$,

$$\begin{aligned} I_{22} &= m_n \sum_{h=k+1}^{r_n} p_{0h} + \mathrm{o}(1) \\ &\leq m_n \sum_{h=k+1}^{r_n} P(|\mathbf{X}_h| > \epsilon a_m, |\mathbf{X}_0| > \epsilon a_m) + \mathrm{o}(1). \end{aligned}$$

We conclude from (3.3) that $\lim_{k\to\infty} \limsup_{n\to\infty} I_{22} = 0$. Finally, since $(\mathbf{X}_t)$ is $\alpha$-mixing and condition (3.2) holds,

$$\lim_{n\to\infty} |I_{23}| \leq \lim_{n\to\infty} m_n \sum_{h=r_n+1}^{\infty} \alpha_h = 0,$$

which completes the proof of the theorem. □

### 3.2. A central limit theorem for $\widehat{P}_m(C)$ and the empirical extremogram

The following central limit theorem is the main result of this paper.

**Theorem 3.2.** *Assume that the conditions of Theorem 3.1 hold with $k_n = n/m_n$, $m_n$ and $r_n$ satisfying $k_n \alpha_{r_n} \to 0$, and $m_n = \mathrm{o}(n^{1/3})$. Then, the central limit theorem*

$$S_n = \left(\frac{n}{m_n}\right)^{1/2} [\widehat{P}_m(C) - m_n P(a_m^{-1}\mathbf{X} \in C)]$$
$$= \left(\frac{m_n}{n}\right)^{1/2} \sum_{i=1}^{n} (I_{\{\mathbf{X}_t/a_m \in C\}} - P(a_m^{-1}\mathbf{X} \in C)) \stackrel{d}{\to} N(0, \sigma^2(C))$$

*holds, where*

$$\sigma^2(C) = \mu(C) + 2 \sum_{h=1}^{\infty} \tau_h(C). \tag{3.9}$$

*The condition $m_n = \mathrm{o}(n^{1/3})$ can be replaced by the condition*

$$\frac{m_n^4}{n} \sum_{j=r_n}^{m_n} \alpha_j \to 0 \quad \textit{and} \quad \frac{m_n r_n^3}{n} \to 0, \tag{3.10}$$

*which is often much weaker.*

The proof of the theorem is given in Section 6. It is based on a standard big-block/small-block argument. Proofs in a similar vein in an extreme value theory context can be found in the literature; see, for example, Rootzén *et al.* [35]. In Section 6, we also propose an estimator of the asymptotic variance $\sigma^2(C)$.

For many examples considered in financial time series and elsewhere, the $\alpha$-mixing rate function $\alpha_j$ decays at an exponential rate. In these cases, one can take $m_n \sim n^{1/2-\delta}$ for some small $\delta > 0$ and $r_n \sim n^{1/8}$. The choice $r_n \sim c \log n$ for some $c > 0$ also fulfills (3.10).

A slight adaptation of the proofs given in Theorems 3.1 and 3.2, in combination with the central limit theorem in Utev [38], shows that these results hold if condition (M) is replaced by the assumption that the process is $\phi$-mixing with a summable rate function $(\phi_t)$.

A related paper on the pre-asymptotic behavior of the empirical extremogram in the case $A = (x, \infty)$ and $B = (y, \infty)$ is Hill [20]; see, in particular, his Theorem 5.4. In contrast to the present paper, he does "not require a model for the bivariate joint tail nor any assumptions concerning the joint tail" (his Remark 15).

### 3.3. Extremogram estimation

In order to derive the limiting distribution of the extremogram estimator, we first consider the large-sample behavior of the ratio estimator of

$$R(C,D) := \frac{\mu(D)}{\mu(C)}$$

given by

$$\widehat{R}_m(C,D) = \frac{\widehat{P}_m(D)}{\widehat{P}_m(C)},$$

where $C$ and $D$ are sets of the type described in Theorem 3.1 with $\mu(C) > 0$. Under the conditions of this theorem, $\widehat{R}_m(C,D)$ is a consistent estimator of $R(C,D)$. In what follows, we study the central limit theorem for this ratio estimator.

Observe that

$$\begin{aligned} &\widehat{R}_m(C,D) - R(C,D) \\ &\quad = \frac{\widehat{P}_m(D)\mu(C) - \widehat{P}_m(C)\mu(D)}{\mu(C)\widehat{P}_m(C)} \\ &\quad = \frac{1 + o_P(1)}{[\mu(C)]^2}[(\mu(C)(\widehat{P}_m(D) - E\widehat{P}_m(D)) - \mu(D)(\widehat{P}_m(C) - E\widehat{P}_m(C))) \\ &\qquad\qquad + (\mu(C)E\widehat{P}_m(D) - \mu(D)E\widehat{P}_m(C))]. \end{aligned} \tag{3.11}$$

The decomposition (3.11) indicates how we have to proceed. First, we must prove a central limit theorem for the first term on the right-hand side. This problem is similar to Theorem 3.2 and requires proving a joint central limit theorem for $(\widehat{P}_m(C), \widehat{P}_m(D))$. For the second term in (3.11), we have, by (3.5),

$$\begin{aligned} &\mu(C)E\widehat{P}_m(D) - \mu(D)E\widehat{P}_m(C) \\ &\quad = m_n[\mu(D)P(a_m^{-1}\mathbf{X} \in C) - \mu(C)P(a_m^{-1}\mathbf{X} \in D)] = o(1). \end{aligned} \tag{3.12}$$

However, for a central limit theorem for $(n/m_n)^{1/2}(\widehat{R}_m(C,D) - R(C,D))$, one needs to know the rate of convergence of (3.12) to zero. This is, in general, a difficult problem which can sometimes be solved when one deals with a specific time series model; see, for example, Section 4.2 in the case of a stochastic volatility model. Alternatively, one could assume conditions on the rate of convergence in the relations $P_m(D) \to \mu(D)$ and $P_m(C) \to \mu(C)$. Such conditions are common in extreme value theory.

We formulate the central limit theorem for the finite-dimensional distributions of the ratio estimator in the following corollary.

**Corollary 3.3.** *Assume the conditions of Theorems 3.1 and 3.2 hold for the sets*[1] $D_1,\dots,D_h,C$ *and the sequence* $(\mathbf{X}_t)$. *Moreover, let* $\mu(C)>0$. *Then,*

$$\left(\frac{n}{m_n}\right)^{1/2}[\widehat{R}_m(C,D_i)-R_m(C,D_i)]_{i=1,\dots,h}\xrightarrow{d}N(\mathbf{0},(\mu(C))^{-4}\mathbf{F}'\mathbf{\Sigma}\mathbf{F}),$$

*where*

$$R_m(C,D_i)=\frac{P(a_m^{-1}\mathbf{X}\in D_i)}{P(a_m^{-1}\mathbf{X}\in C)}$$

*and where* $\mathbf{\Sigma}$ *and* $\mathbf{F}$ *are defined in (3.15) and (3.18), respectively. If, in addition,*

$$\lim_{n\to\infty}\sqrt{nm_n}[\mu(D_i)P(a_m^{-1}\mathbf{X}\in C)-\mu(C)P(a_m^{-1}\mathbf{X}\in D_i)]=0,\qquad i=1,\dots,h, \tag{3.13}$$

*then*

$$\left(\frac{n}{m_n}\right)^{1/2}[\widehat{R}_m(C,D_i)-R(C,D_i)]_{i=1,\dots,h}\xrightarrow{d}N(\mathbf{0},(\mu(C))^{-4}\mathbf{F}'\mathbf{\Sigma}\mathbf{F}).$$

**Proof.** In order to ease notation, we set $D_{h+1}=C$. We show the central limit theorem

$$\mathbf{S}_n=\left(\frac{m_n}{n}\right)^{1/2}\sum_{t=1}^{n}\begin{pmatrix} I_{D_1}(\mathbf{X}_t/a_m)-P(a_m^{-1}\mathbf{X}\in D_1)\\ \vdots \\ I_{D_{h+1}}(\mathbf{X}_t/a_m)-P(a_m^{-1}\mathbf{X}\in D_{h+1})\end{pmatrix}\xrightarrow{d}N(\mathbf{0},\mathbf{\Sigma}), \tag{3.14}$$

where

$$\mathbf{\Sigma}=\begin{pmatrix} \sigma^2(D_1) & r_{D_1,D_2} & r_{D_1,D_3} & \cdots & r_{D_1,D_{h+1}}\\ \vdots & \vdots & \vdots & \ddots & \vdots\\ r_{D_1,D_{h+1}} & r_{D_2,D_{h+1}} & r_{D_3,D_{h+1}} & \cdots & \sigma^2(D_{h+1})\end{pmatrix} \tag{3.15}$$

and, for $i\neq j$,

$$r_{D_i,D_j}=\mu(D_i\cap D_j)+2\sum_{h=1}^{\infty}\mu_{h+1}(D_i\times(\overline{\mathbb{R}}_0^d)^{h-2}\times D_j). \tag{3.16}$$

By the Cramér–Wold device, it suffices to show the central limit theorem for any linear combination

$$\mathbf{z}'\mathbf{S}_n\xrightarrow{d}N(\mathbf{0},\mathbf{z}'\mathbf{\Sigma}\mathbf{z}),\qquad \mathbf{z}\in\mathbb{R}^{h+1}.$$

[1]This means that the set $C$ in Theorems 3.1 and 3.2 has to be replaced by the $D_i$'s.

The same ideas as in the proof of Theorem 3.2 show that it suffices to prove the central limit theorem for $k_n$ i.i.d. copies of

$$T_n(\mathbf{z}) = (m_n/n)^{1/2} \sum_{t=1}^{m_n} \sum_{i=1}^{h+1} z_i (I_{D_i}(\mathbf{X}_t/a_m) - P(a_m^{-1}\mathbf{X} \in D_i)).$$

By the central limit theorem for triangular arrays, one needs to verify that for every $\epsilon > 0$,

$$k_n E(T_n^2(\mathbf{z}) I_{\{|T_n(\mathbf{z})|>\epsilon\}}) \to 0. \tag{3.17}$$

This follows from Markov's inequality and (6.4) when $m_n = \mathrm{o}(n^{1/3})$:

$$k_n E T_n^2(\mathbf{z}) I_{\{|T_n(\mathbf{z})|>\epsilon\}} \le c m_n^2 P(|T_n(\mathbf{z})| > \epsilon) \le c m_n^2 k_n^{-1} = \mathrm{o}(1).$$

If the conditions (3.10) are met, then the argument at the end of the proof given in Section 6 can be used to establish (3.17). This proves the central limit theorem (3.14).

We observe that

$$\begin{aligned}&\widehat{R}_m(D_{h+1}, D_i) - R_m(D_{h+1}, D_i)\\ &\quad = \frac{1 + \mathrm{o}_P(1)}{[\mu(D_{h+1})]^2}((\widehat{P}_m(D_i) - P_m(D_i))\mu(D_{h+1}) - (\widehat{P}_m(D_{h+1}) - P_m(D_{h+1}))\mu(D_i)).\end{aligned}$$

Hence,

$$\begin{aligned}\left(\frac{n}{m_n}\right)^{1/2} [\widehat{R}_m(D_{h+1}, D_1) - R_m(D_{h+1}, D_1)]_{i=1,\dots,h} &= \frac{1 + \mathrm{o}_P(1)}{(\mu(D_{h+1}))^2} \mathbf{F}\mathbf{S}_n\\ &\stackrel{d}{\to} N(\mathbf{0}, (\mu(D_{h+1}))^{-4}\mathbf{F}'\mathbf{\Sigma}\mathbf{F}),\end{aligned}$$

where

$$\mathbf{F} = \begin{pmatrix} \mu(D_{h+1}) & 0 & 0 & \cdots & 0 & -\mu(D_1) \\ 0 & \mu(D_{h+1}) & 0 & \cdots & 0 & -\mu(D_2) \\ \vdots & \vdots & \vdots & \ddots & \vdots & \vdots \\ 0 & 0 & 0 & \cdots & \mu(D_{h+1}) & -\mu(D_h) \end{pmatrix}. \tag{3.18}$$

This proves the result. □

Recall that, for subsets $A$ and $B$ of $\overline{\mathbb{R}}_0^d$ that are bounded away from zero and $\mu(\partial A) = \mu(\partial B) = 0$, $\mu(A) > 0$, the extremogram at lag $h$ is defined by

$$\rho_{AB}(h) = \lim_{n\to\infty} \frac{P(a_n^{-1}\mathbf{X}_0 \in A, a_n^{-1}\mathbf{X}_h \in B)}{P(a_n^{-1}\mathbf{X}_0 \in A)}.$$

A natural estimator of $\rho_{AB}$ is the *empirical extremogram* defined by

$$\hat{\rho}_{AB}(i) = \frac{\sum_{t=1}^{n-i} I_{\{a_m^{-1}\mathbf{X}_t \in A, a_m^{-1}\mathbf{X}_{t+i} \in B\}}}{\sum_{t=1}^{n} I_{\{a_m^{-1}\mathbf{X}_t \in A\}}}, \qquad i = 0, 1, \ldots.$$

This estimate can be recast as a ratio estimator by introducing the vector process

$$\mathbf{Y}_t = \operatorname{vec}(\mathbf{X}_t, \ldots, \mathbf{X}_{t+h})$$

consisting of stacking $h+1$ consecutive values of the time series $(\mathbf{X}_t)$. Now, the sets $C$ and $D_0, \ldots, D_h$ specified in Corollary 3.3 are defined via the relations $C = A \times \overline{\mathbb{R}}_0^{dh}$, $D_0 = A \cap B \times \overline{\mathbb{R}}_0^{dh}$ and $D_i = A \times \overline{\mathbb{R}}_0^{d(i-1)} \times B \times \overline{\mathbb{R}}_0^{d(h-i)}$ for $i \geq 1$. With this conversion, Corollary 3.3 can be applied to the $(\mathbf{Y}_t)$ sequence directly to show that $\hat{\rho}_{AB}(i)$, centered by the *pre-asymptotic value* of the extremogram defined by

$$\rho_{AB,m}(i) = \frac{P(a_m^{-1}\mathbf{X}_0 \in A, a_m^{-1}\mathbf{X}_i \in B)}{P(a_m^{-1}\mathbf{X} \in A)}, \tag{3.19}$$

is asymptotically normal. On the other hand, if the bias condition (3.13) is met, then one can center the empirical extremogram by its true value and still retain the asymptotic normality. For completeness, we record these results as the following corollary.

**Corollary 3.4.** *Assume that the conditions of Corollary 3.3 are satisfied for the sequence $(\mathbf{Y}_t)$. Then,*

$$\left(\frac{n}{m_n}\right)^{1/2} [\hat{\rho}_{AB}(i) - \rho_{AB,m}(i)]_{i=0,1,\ldots,h} \xrightarrow{d} N(\mathbf{0}, (\mu(A))^{-4}\mathbf{F}'\mathbf{\Sigma}\mathbf{F}). \tag{3.20}$$

*Moreover, if (3.13) is satisfied, then (3.20) holds with $\rho_{AB,m}(i)$ replaced by $\rho_{AB}(i)$.*

## 3.4. An empirical example

For illustrative purposes, we compute the extremogram for the high-frequency financial time series consisting of 35135 five minute returns of the foreign exchange rate USD-DEM. The data, which was provided in processed form by Olsen and Associates (Zürich) at their Second Conference on High Frequency Financial Data in 1995, was one of the first widely disseminated high frequency data sets. Aside from the choice of $A$ and $B$, the most problematic issue in computing the extremogram is the selection of a suitable threshold $a_m$. The choice of threshold is always a thorny issue in extreme value theory, from estimating the tail index of regular variation to fitting a generalized Pareto distribution. We tried several different choices of threshold based on various empirical quantiles of the absolute values of the data. Of course, there is the typical bias-variance trade-off in this selection with a large $m$ corresponding to a smaller bias, but larger variance, and vice versa for a moderate value of $m$. In the upper-left panel of Figure 1, we plot the sample extremogram for lags 1–100 using $A = B = (-\infty, -1) \cup (1, \infty)$ with several choices of threshold $a_m$ starting with the 0.97 empirical quantile of the absolute values (solid line) and continuing with the 0.98, 0.99, 0.995 and 0.999 quantiles (dotted and dashed lines). With the exception of the 0.995 and 0.999 quantiles, the extremogram has roughly the same value for all lags. The extremogram based on the 0.995 quantile is also similar, at least for the first 38 lags, and drops off to zero due to a lack of pairs of lagged observations that exceed this large quantile. In view of this robustness of the plots for this range of thresholds, we will choose $a_m$ to be the 0.98 quantile in all of the remaining extremogram plots.

The extremogram of the returns with $A = B = (1, \infty)$ is displayed in the upper-right panel of Figure 1. Note that the extremogram alternates between large values at even lags and smaller values at odd lags. Like many high frequency data sets, the autocorrelation function for this time series alternates between positive and negative values. To further investigate this alternating behavior of the extremogram, which may be due, in part, to an artifact of the processing of the data by Olsen, we fitted an AR model to the data. The best fitting AR model, based on minimizing the AICC, is of order 18. We then refined this model by selecting the best subset model which ended up having significant non-zero coefficients at lags 1, 2, 3, 5, 6, 7, 11, 13, 14, 16 and 18. The lag-1 coefficient, which was much larger than the other lags, was –0.6465. So, this alternating character of the extremogram is consistent with the extremogram of an AR(1) process with negative coefficient described in Section 2.5. In the middle-left panel, we plot the extremogram of the residuals from the subset AR(18) model fit with $A = B = (1, \infty)$. As a baseline, we have also plotted the horizontal line that one would expect for the extremogram if the data were in fact independent and the threshold was the 0.98 quantile of the absolute values. Note that the values are now significantly smaller than those for the returns. Some extremal dependence still remains in the residuals, at least for small lags. This behavior is an indication of the presence of nonlinearity in the data. In fact, the ACFs of the absolute values and squares of the residuals are highly significant. We were moderately successful in eliminating part of this nonlinearity by fitting a GARCH$(1, 1)$ model with $t$-noise to the residuals. The extremogram of the GARCH residuals, which still exhibits some ACF in the absolute values but none in the squares, is displayed in the middle-right panel of Figure 1. Based on this extremogram, there is little remaining extremal dependence in the GARCH residuals. As a last check on this modeling exercise, we simulated a realization of the time series based on the fitted model. In other words, we generated a time series from the GARCH model and then passed it through the fitted AR filter. The extremogram of this simulated series is displayed at the bottom of Figure 1. It displays similar features to the original extremogram (top-right panel of Figure 1), but the dependence is not quite as strong or persistent as for the original data.

In Figure 2, we chose $A = (1, \infty)$ and $B = (-\infty, -1)$ for computing the extremograms for the return data (left) and the residuals from the AR(18) fit (right). For the return data, the extremogram at the first lag has a large positive value and alternates at the odd and even lags. On the other hand, for the residuals, there is a real difference in shape of the extremogram from that displayed in Figure 1.

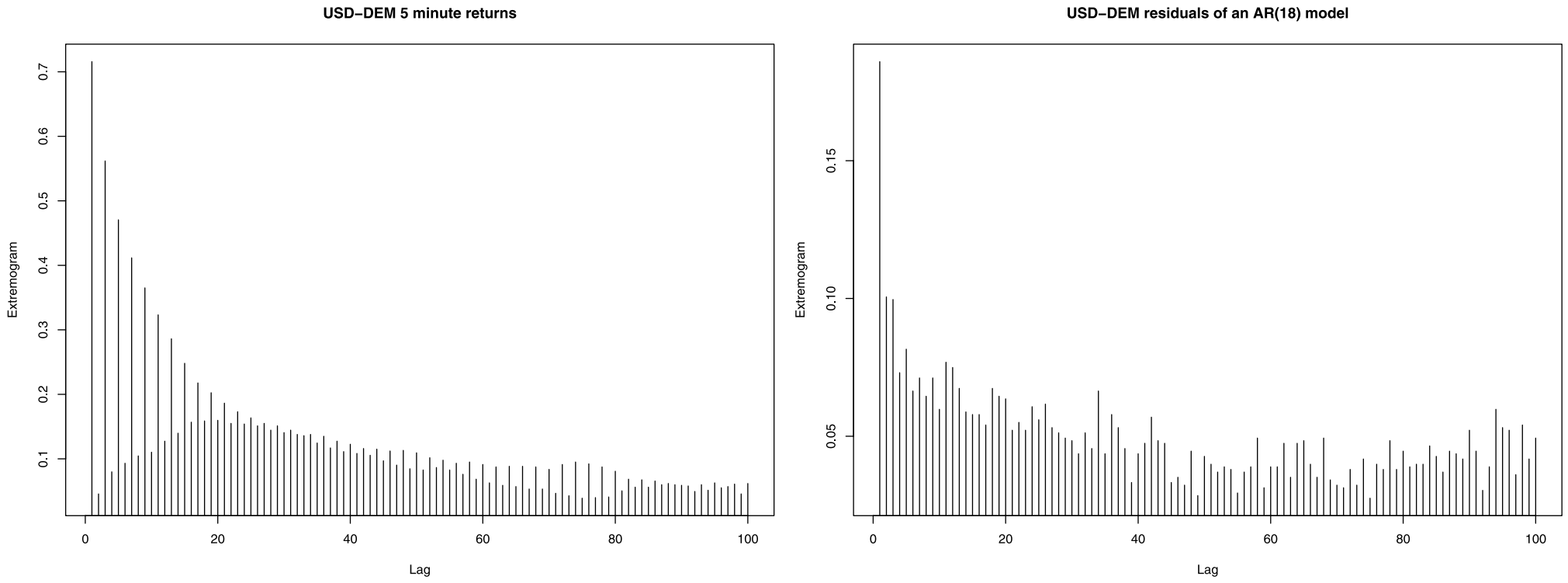

**Figure 2.** The empirical extremogram (left) of a sample of five-minute return data of the foreign exchange rate USD-DEM and the residuals (right) from a fitted AR(18) model with $A=(1,\infty)$, $B=(-\infty,-1)$; see Section 3.4 for more information.

# 4. Application to GARCH and SV models

## 4.1. The GARCH process

Assume the conditions of Section 2.3 hold for a regularly varying, strictly stationary GARCH$(1,1)$ process with index $\alpha>0$. The GARCH process is $\beta$-mixing, hence $\alpha$-mixing, with geometric rate under general conditions on the noise; see Boussama [5] and Mokkadem [27]. In the GARCH$(1,1)$ case, these conditions hold provided the density of the noise variables $Z_t$ is positive in some neighborhood of the origin.

In what follows, we assume that $(X_t)$ is $\alpha$-mixing with a geometric rate function $\alpha_t \le ca^t$ for some $a\in(0,1)$, $c>0$. First, recall that the $\kappa$ chosen in Section 2.3 satisfies $\kappa\in(0,\alpha/2)$ and $EC^\kappa<1$, where $C=\alpha_1 Z_1^2+\beta_1$. Second, the normalizing constants $a_n$ are chosen such that $P(|X|>a_n)\sim n^{-1}$. In particular, $a_n\sim cn^{1/\alpha}$. Now, select $m_n=n^\delta$ ($\delta\in(0,1)$) and $r_n=n^\gamma$ for $\gamma<\min((1-\delta)/3,\delta(2\kappa/\alpha))$. With these choices of $m_n,r_n$, we verify conditions (M) and (3.10). The mixing condition (3.2) is straightforward to check since $(X_t)$ is $\alpha$-mixing with geometric rate. To check (3.3), it suffices to show (3.4). Using the recursion in (2.2), we have

$$\sigma_t^2 = C_t\cdots C_1\sigma_0^2+\alpha_0\sum_{i=1}^t C_t\cdots C_{i+1},\qquad t\ge 0. \tag{4.1}$$

If we assume, for the sake of simplicity, that $\epsilon=1$, then (4.1) and Markov's inequality imply that

$$\begin{aligned}
&P(|X_h|>a_m \mid |X_0|>a_m)\\
&\quad\le P(Z_h^2C_h\cdots C_1\sigma_0^2>a_m^2/2 \mid |X_0|>a_m)+P\left(Z_h^2\alpha_0\sum_{i=1}^h C_h\cdots C_{i+1}>a_m^2/2\right)
\end{aligned}$$

$$\leq \frac{P(Z_h^2 C_h \cdots C_1 \sigma_0^2 > a_m^2/2)}{P(|X_0| > a_m)} + c(a_m^2/2)^{-\kappa} E|Z|^{2\kappa} \sum_{i=1}^{h} (EC^\kappa)^{h-i}$$
$$= I_1(h) + I_2(h).$$

With the prescribed choices of $r_n$ and $m_n$, we have

$$\sum_{h=k+1}^{r_n} I_2(h) \leq c r_n a_m^{-2\kappa} \to 0, \qquad n \to \infty.$$

By again applying Markov's inequality and Karamata's theorem (see Bingham *et al.* [4]),

$$\begin{aligned}\sum_{h=k+1}^{r_n} I_1(h) &= \sum_{h=k+1}^{r_n} \frac{P(C_h \cdots C_1 X_0^2 > a_m^2/2)}{P(|X| > a_m)} \\ &\leq c \frac{E|X|^\kappa I_{\{|X|>a_m\}}}{a_m^\kappa P(|X| > a_m)} \sum_{h=k+1}^{r_n} (EC^\kappa)^h \\ &\leq cE(C^\kappa)^k.\end{aligned}$$

Appropriately combining the above facts shows that (3.3) is satisfied.

Finally, it is straightforward to check (3.10) from the choices of $\delta$ and $\gamma$. Hence, the conclusions of Theorems 3.1 and 3.2 apply to the GARCH$(1,1)$ process. To date, we have not been able to verify the bias condition (3.13) of Corollary 3.4. One needs a more precise estimate than is currently known for the tail distribution of $\sigma_t^2$; see Goldie [18] for results in this direction.

## 4.2. The stochastic volatility process

The stochastic volatility process $(X_t)$ has regularly varying finite-dimensional distributions with index $\alpha$ if the multiplicative noise $(Z_t)$ is regularly varying with index $\alpha$ and the volatility $\sigma_t$ has a moment of order $\alpha + \epsilon$, $\epsilon > 0$. In fact, for the following argument, we will assume that $\sigma_t$ has a finite $4\alpha$th moment. In addition, we assume that the mixing condition (3.2) is satisfied for $(\sigma_t)$. If the sequence $(\sigma_t)$ is $\alpha$-mixing with rate function $(\alpha_h)$, then $(X_t)$ has rate function $(4\alpha_h)$, and hence $(X_t)$ also satisfies (3.2).

We next dispense with condition (3.3) with $\epsilon = 1$; the general case $\epsilon > 0$ is completely analogous. Using the independence of $(Z_t)$ and $(\sigma_t)$, an application of Markov's inequality for $p < \alpha$ yields

$$\begin{aligned}&m_n \sum_{h=k}^{r_n} P(a_m^{-1}|\mathbf{X}_0| > 1, a_m^{-1}|\mathbf{X}_h| > 1) \\ &\quad \leq m_n \sum_{h=k}^{r_n} P(\max(\sigma_0, \sigma_1, \sigma_h, \sigma_{h+1}) \min(|Z_0|, |Z_1|, |Z_h|, |Z_{h+1}|) > a_m)\end{aligned}$$

$$= m_n \sum_{h=k}^{r_n} E[P(\max(\sigma_0, \sigma_1, \sigma_h, \sigma_{h+1})Z > a_m | (\sigma_t))]^4$$
$$\le c m_n r_n (E\sigma^p)^4 a_m^{-4p}.$$

Since $r_n = o(m_n)$, the right-hand side vanishes if $p$ is chosen close to $\alpha$.

We now turn to the problem of verifying the bias condition (3.13) so that we can apply the limit theory of Corollaries 3.3 and 3.4 to the empirical extremogram. To this end, we assume, for convenience, that $\log \sigma_t$ is a stationary Gaussian sequence with mean zero and unit variance. Choose $m_n = n^\gamma$ for some $\gamma \in (1/3, 1)$ and suppose the mixing function decays sufficiently fast so that (3.10) holds. For example, if $\alpha_t$ decays geometrically, one could take $r_n = (1-\gamma)/4$. If we choose the sets $A = (1, \infty) \times (0, \infty)^{h-1}$ and $B = (1, \infty) \times (0, \infty)^{h-2} \times (1, \infty)$, then $\mu_h(A) > 0$ and $\mu_h(B) = 0$. Set $s_n = n^\delta$ for some $0 < \delta < (3\gamma - 1)/(4\alpha)$ and note that

$$P(a_m^{-1} X_0 > 1, a_m^{-1} X_h > 1) \le P(\max(\sigma_0, \sigma_h) \min(Z_0, Z_1) > a_m)$$
$$\le P(\max(\sigma_0, \sigma_h) > s_n) + P(\min(Z_0, Z_1) > a_m / s_n).$$

Since $(m_n n)^{1/2} = n^{0.5(1+\gamma)}$, an application of Markov's inequality yields that for any $k > 0$,

$$(m_n n)^{1/2} P(\max(\sigma_0, \sigma_h) > s_n) \le n^{0.5(1+\gamma)} 2E(\sigma^k) s_n^{-k}.$$

The right-hand side converges to zero for $k$ sufficiently large. On the other hand, for any $\epsilon > 0$,

$$(m_n n)^{1/2} [P(Z > a_m / s_n)]^2 \le n^{0.5(1+\gamma)} (n^{\gamma/\alpha} / n^\delta)^{-2\alpha+\epsilon}.$$

The right-hand side converges to zero for small $\epsilon$. This shows that (3.13) is satisfied for a stochastic volatility model and sets $A, B$ as specified. Applying Corollary 3.4, we conclude that

$$(n/m_n)^{1/2} [\hat{\rho}_{AB}(1) - \rho_{AB}(1)] \stackrel{d}{\to} N(0, \sigma^2(B) \mu^{-2}(A)),$$

where, of course, $\rho_{AB}(1) = 0$ in this case and $\sigma^2(B) = 0$. Therefore, we get a degenerate limit for this choice of $A$ and $B$.

As a second example, let $A = (1, \infty) \times (0, \infty)$ and $B = \{(x_1, x_2)\colon\ L < x_1 - x_2 < U, x_1 \ge 0, x_2 \ge 0\}$, where $L < 1 < U$. With $\mathbf{X}_t = (X_t, X_{t+1})$, a straightforward calculation shows that $\mu(A) = 1$ and $\mu(B) = L^{-\alpha} - U^{-\alpha}$. Since the measure $\mu$ only concentrates on the two coordinate axes, $A \cap B$ intersects the $x_2$-axis in the empty set and intersects the $x_1$-axis in the interval $[1, U]$. Hence, $\gamma_{AB}(0) = 1 - U^{-\alpha}$. Since the limiting measure of $(X_0, X_1, X_h, X_{h+1})$ concentrates on the four coordinate axes,

$$nP(a_n^{-1} \mathbf{X}_0 \in A, a_n^{-1} \mathbf{X}_h \in B) \to \gamma_{AB}(h) = 0$$

for $h \geq 1$. On the other hand, $\{(x_1, x_2, x_3)\colon\ (x_1, x_2) \in B \text{ and } (x_2, x_3) \in B\}$ has empty intersection with the three coordinate axes and hence $\gamma_{BB}(1) = 0$. More generally, we have

$$\gamma_{BB}(0) = L^{-\alpha} - U^{-\alpha}, \qquad \gamma_{BB}(h) = \gamma_{AA}(h) = 0 \qquad \text{for } h > 0.$$

Using Corollary 3.3, we have

$$(n/m_n)^{1/2}(\hat{\rho}_{AB}(0) - \rho_{AB,m}(0)) \xrightarrow{d} N(0, 1 - U^{-\alpha}).$$

## 5. Some spectral analysis

The extremogram $\gamma_{CC} = (\tau_h(C))$ with $\tau_h(C) = \tau_{-h}$ and $\tau_0(C) = \mu(C)$ as defined in Theorem 3.1 is an asymptotic covariance function. If it is summable, then the function

$$f(\lambda) = \sum_{h \in \mathbb{Z}} \tau_h(C) \mathrm{e}^{ih\lambda} = \mu(C) + 2 \sum_{h=1}^{\infty} \cos(\lambda h) \tau_h(C), \qquad \lambda \in [0, \pi],$$

defines the corresponding spectral density which, in turn, determines $\gamma_{CC}$. The sample version of the spectral density $f$ is given by the periodogram

$$I_{nC}(\lambda) = \frac{m_n}{n} \left| \sum_{t=1}^{n} \widetilde{I}_t \mathrm{e}^{it\lambda} \right|^2 = \widehat{\gamma}_n(0) + 2 \sum_{h=1}^{n-1} \cos(\lambda h) \widehat{\gamma}_n(h), \qquad \lambda \in [0, \pi],$$

where

$$p_0 = P(\mathbf{X}_t / a_m \in C), \qquad I_t = I_{\{\mathbf{X}_t / a_m \in C\}} \quad \text{and} \quad \widetilde{I}_t = I_t - EI_t = I_t - p_0,$$

and

$$\widehat{\gamma}_n(h) = \frac{m_n}{n} \sum_{t=1}^{n-h} \widetilde{I}_t \widetilde{I}_{t+h}, \qquad \widetilde{\gamma}_n(h) = \frac{m_n}{n} \sum_{t=1}^{n-h} I_t I_{t+h}, \qquad h \geq 0,$$

are analogs of the sample autocovariance function of a stationary sequence. Since it is common to evaluate the quantities $I_{nC}(\lambda)$ at the Fourier frequencies $\lambda = \lambda_k = 2\pi k/n \in (0, \pi)$ and $\sum_{t=1}^{n} \mathrm{e}^{it\lambda_k} = 0$, one can define $I_{nC}(\lambda_k)$ with the $\widetilde{I}_t$'s replaced by the $I_t$'s which do not contain the unknown probability $p_0$. However, for the calculations which involve mixing conditions, it is crucial to use the given definition of $I_{nC}(\lambda)$ with the centered quantities $\widetilde{I}_t$.

In what follows, we will mostly deal with the *lag-window estimator* (see Brockwell and Davis [8]) or the *truncated periodogram* $\widehat{f}_{nC}(\lambda)$ defined by

$$\widehat{f}_{nC}(\lambda) = \widehat{\gamma}_n(0) + 2 \sum_{h=1}^{r_n} \cos(\lambda h) \widetilde{\gamma}_n(h), \qquad \lambda \in [0, \pi],$$

where $r_n \to \infty$ and $r_n/m_n \to 0$ as $n \to \infty$ has the same interpretation as in the previous sections. Truncated estimators of the form $\widehat{f}_{nC}$ are commonly used in the spectral analysis of stationary time series; see, for example, Brockwell and Davis [8] and Priestley [29]. A major reason for this is that, unlike the periodogram of a stationary time series, the truncated periodogram $\widehat{f}_{nC}$ is a consistent estimator of the spectral density. In our setting, we show below that $\widehat{f}_{nC}$ remains a consistent estimator of $f(\lambda)$. Based on background calculations, it appears that $I_{nC}(\lambda)$ is not consistent.

**Theorem 5.1.** *Assume the mixing condition* (M) *for the regularly varying, strictly stationary sequence* $(\mathbf{X}_t)$ *with index* $\alpha > 0$ *and that the products* $C^k \subset \overline{\mathbb{R}}_0^{dk}$ *are continuity sets with respect to the limiting measures* $\mu_k$, $k = 1, 2, \ldots,$ *occurring in the definition of regular variation. Then,*

$$\lim_{n\to\infty} EI_{nC}(\lambda) = \lim_{n\to\infty} \widehat{f}_{nC}(\lambda) = f(\lambda), \qquad \lambda \in (0, \pi). \tag{5.1}$$

*In addition, if* $m_n r_n^2 = \mathrm{O}(n)$, *then we also have*

$$\lim_{n\to\infty} E[(\widehat{f}_{nC}(\lambda) - f(\lambda))^2] = 0, \qquad \lambda \in (0, \pi). \tag{5.2}$$

*This means that the estimator* $\widehat{f}_{nC}(\lambda)$ *of the spectral density* $f(\lambda)$ *is asymptotically unbiased and mean-square consistent.*

The rate of convergence in (5.2) cannot be derived unless one assumes conditions similar to (3.13).

**Proof.** In what follows, it will be convenient to use the notation

$$p_0 = P(\mathbf{X}/a_m \in C), \qquad p_{st} = P(\mathbf{X}_s/a_m \in C, \mathbf{X}_t/a_m \in C), \qquad \ldots,$$
$$\widetilde{I}_{st} = \widetilde{I}_s \widetilde{I}_t, \qquad \ldots, \qquad \widetilde{p}_{st} = E\widetilde{I}_{st}, \qquad \ldots.$$

We will exploit the following auxiliary result. The proof is completely analogous to the proof of Theorem 3.1 and is therefore omitted.

**Lemma 5.2.** *Under the conditions of Theorem 5.1,*

$$E\widehat{\gamma}_n(h) \sim \widetilde{\gamma}_n(h) \to \tau_h(C),$$
$$\operatorname{var}(\widehat{\gamma}_n(h)) \sim \operatorname{var}(\widetilde{\gamma}_n(h)) \sim \frac{m_n}{n}\left[\tau_h(C) + 2\sum_{t=1}^{\infty} \tau_{0ht,t+h}(C)\right], \qquad h \geq 0,$$

*where, for* $h \geq 0$ *and* $t \geq 0$, $\tau_{0ht,t+h}(C) = \lim_{n\to\infty} m p_{0ht,t+h}$.

We start by considering the expectation of the periodogram. We have, for fixed $k\ge 1$,

$$\begin{aligned}
EI_{nC}(\lambda) &= m_n p_0 + 2m_n\sum_{h=1}^{n-1}(1-n^{-1}h)\cos(\lambda h)\widetilde{p}_{0h}\\
&= m_n p_0 + 2m_n\sum_{h=1}^{n-1}(1-n^{-1}h)\cos(\lambda h)[p_{0h}-p_0^2]\\
&= \mu(C) + 2\sum_{h=1}^{k}\cos(\lambda h)\tau_h(C)\\
&\quad + 2m_n\left(\sum_{h=k+1}^{r_n}+\sum_{h=r_n+1}^{n-1}\right)(1-n^{-1}h)\cos(h\lambda)[p_{0h}-p_0^2] + \mathrm{o}(1).
\end{aligned}$$

Using condition (M),

$$\lim_{k\to\infty}\limsup_{n\to\infty} m_n\left|\sum_{h=k+1}^{r_n}(1-n^{-1}h)\cos(h\lambda)[p_{0h}-p_0^2]\right| \le \lim_{k\to\infty}\limsup_{n\to\infty} m_n\sum_{h=k+1}^{r_n}p_{0h}=0,$$

$$\lim_{n\to\infty} m_n\left|\sum_{h=r_n+1}^{n-1}(1-n^{-1}h)\cos(h\lambda)[p_{0h}-p_0^2]\right| \le \limsup_{n\to\infty} m_n\sum_{h=r_n+1}^{\infty}\alpha_h=0.$$

The relation $\lim_{n\to\infty}E\widehat{f}_{nC}(\lambda)=f(\lambda)$ is derived in the same way.

We conclude from Lemma 5.2 that for any $k\ge 1$,

$$\mathrm{var}\left((\widetilde{\gamma}_n(0)-\mu(C)) + 2\sum_{h=1}^{k}\cos(\lambda h)(\widetilde{\gamma}_n(h)-\tau_h(C))\right)\to 0.$$

Therefore, it suffices for (5.2) to show that

$$\lim_{k\to\infty}\limsup_{n\to\infty}\mathrm{var}\left(\sum_{h=k}^{r_n}\cos(\lambda h)\widetilde{\gamma}_n(h)\right)=0.$$

It suffices to bound the expression

$$I=\sum_{h=k}^{r_n}\sum_{l=k}^{r_n}|\operatorname{cov}(\widetilde{\gamma}_n(h),\widetilde{\gamma}_n(l))| \le \frac{m_n^2}{n^2}\sum_{h=k}^{r_n}\sum_{l=k}^{r_n}\sum_{t=1}^{n-h}\sum_{s=1}^{n-l}|\operatorname{cov}(I_{t,t+h},I_{s,s+l})|.$$

We have

$$I\le 2\frac{m_n^2}{n}\sum_{h=k}^{r_n}\sum_{l=k}^{r_n}\sum_{r=0}^{n}|\operatorname{cov}(I_{0h},I_{r,r+l})|$$

$$\leq 2\frac{m_n^2}{n}\sum_{h=k}^{r_n}\sum_{l=k}^{r_n}\sum_{r=0}^{2r_n-1}\sqrt{\operatorname{var}(I_{0h})\operatorname{var}(I_{0l})}+2\frac{m_n^2}{n}\sum_{h=k}^{r_n}\sum_{l=k}^{r_n}\sum_{r=2r_n}^{\infty}\alpha_{r-h}$$
$$\leq 4\frac{m_n^2 r_n}{n}\sum_{h=k}^{r_n}\sum_{l=k}^{r_n}[p_{0h}+p_{0l}+2p_0^3]+2\frac{m_n r_n^2}{n}\left[m_n\sum_{r=r_n}^{\infty}\alpha_r\right]$$
$$\leq 8\frac{m_n r_n}{n}\left[m_n\sum_{h=k}^{r_n}p_{0h}\right]+8(m_n p_0)^2 p_0\frac{r_n^3}{n}+2\frac{m_n r_n^2}{n}\left[m_n\sum_{r=m_n-r_n}^{\infty}\alpha_r\right].$$

Condition (M) and the growth restrictions $m_n r_n^2 = \mathrm{O}(n)$ yield the desired result

$$\lim_{k\to\infty}\limsup_{n\to\infty} I = 0.$$

□

# 6. Proof of Theorem 3.2

**Proof.** We use the same notation as in Section 5 and write

$$Y_{nt} = (m_n/n)^{1/2}(I_{\{\mathbf{X}_t/a_m\in C\}} - p_0) = (m_n/n)^{1/2}\widetilde{I}_t, \qquad t = 1,\ldots,n.$$

In order to prove the result, we will use the technique of small/large blocks which is well known in the asymptotic theory for sums of dependent random variables. For simplicity, we will assume that $n/m_n = k_n$ is an integer. The non-integer case does not present any additional difficulties, but requires additional bookkeeping. We introduce the index sets

$$I_{ni} = \{(i-1)m_n + 1,\ldots, im_n\}, \qquad i = 1,\ldots,k_n.$$

By $\widetilde{I}_{ni}$, we denote the index set which consists of all elements of $I_{ni}$ but the first $r_n$ elements and we also write $J_{ni} = I_{ni}\setminus\widetilde{I}_{ni}$. Since $r_n/m_n \to 0$ and $m_n\to\infty$, the sets $\widetilde{I}_{ni}$ are non-empty for large $n$. For any index set $B$ of the integers, we write

$$S_n(B) = \sum_{j\in B} Y_{nj}.$$

We first show that

$$\operatorname{var}\left(\sum_{l=1}^{k_n} S_n(J_{nl})\right)\to 0. \tag{6.1}$$

We have

$$\operatorname{var}\left(\sum_{l=1}^{k_n} S_n(J_{nl})\right)\leq k_n\operatorname{var}(S_n(J_{n1})) + 2k_n\sum_{h=1}^{k_n-1}|\operatorname{cov}(S_n(J_{n1}),S_n(J_{n,h+1}))|$$
$$= P_1 + P_2.$$

We observe that, by (3.3) and since $r_n/m_n \to 0$,

$$P_1 \sim \left[ r_n p_0 + 2 \sum_{h=1}^{r_n-1} (r_n - h) p_{0h} \right] - (r_n p_0)^2$$
$$\leq 2(r_n/m_n) m_n \sum_{h=k+1}^{r_n-1} p_{0h} + \mathrm{o}(1) = \mathrm{o}(1).$$

Moreover, for positive constants $c$,

$$P_2 = \sum_{h=1}^{k_n-1} \left| \sum_{t \in J_{n1}} \sum_{s \in J_{n,h+1}} \operatorname{cov}(I_t, I_s) \right| \leq \sum_{h=1}^{k_n-1} \sum_{t \in J_{n1}} \sum_{s \in J_{n,h+1}} \alpha_{s-t}$$
$$\leq c r_n \sum_{h=m_n-r_n+1}^{\infty} \alpha_h \leq c m_n \sum_{h=r_n+1}^{\infty} \alpha_h = \mathrm{o}(1).$$

This proves (6.1).

Condition (6.1) implies that $S_n$ and $\sum_{i=1}^{k_n} S_n(\widetilde{I}_{ni})$ have the same limit distribution, provided such a limit exists. Let $\widetilde{S}_n(\widetilde{I}_{ni})$, $i = 1, \ldots, k_n$, be i.i.d. copies of $S_n(\widetilde{I}_{n1})$. In what follows, we use a classical idea due to Bernstein dating back to the 1920s. Iterated use of the definition of $\alpha$-mixing and standard results for strong mixing sequences (see Doukhan [14]) yield, for any $t \in \mathbb{R}$,

$$\left| E \prod_{l=1}^{k_n} \mathrm{e}^{itS_n(\widetilde{I}_{nl})} - E \prod_{l=1}^{k_n} \mathrm{e}^{it\widetilde{S}_n(\widetilde{I}_{nl})} \right|$$
$$= \left| \sum_{l=1}^{k_n} E \prod_{s=1}^{l-1} \mathrm{e}^{itS_n(\widetilde{I}_{ns})} \left( \mathrm{e}^{itS_n(\widetilde{I}_{nl})} - \mathrm{e}^{it\widetilde{S}_n(\widetilde{I}_{nl})} \right) \prod_{s=l+1}^{k_n} \mathrm{e}^{it\widetilde{S}_n(\widetilde{I}_{ns})} \right| \tag{6.2}$$
$$\leq \sum_{l=1}^{k_n} \left| E \prod_{s=1}^{l-1} \mathrm{e}^{itS_n(\widetilde{I}_{ns})} \left( \mathrm{e}^{itS_n(\widetilde{I}_{nl})} - \mathrm{e}^{it\widetilde{S}_n(\widetilde{I}_{nl})} \right) \prod_{s=l+1}^{k_n} \mathrm{e}^{it\widetilde{S}_n(\widetilde{I}_{ns})} \right| \leq k_n 4 \alpha_{r_n}.$$

By assumption, the right-hand side converges to zero as $n \to \infty$. Hence, $\sum_{l=1}^{k_n} S_n(\widetilde{I}_{nl})$ and $\sum_{l=1}^{k_n} \widetilde{S}_n(\widetilde{I}_{nl})$ have the same limits in distribution (provided these limits exist). Let $\widetilde{S}_n(I_{ni})$, $i = 1, \ldots, k_n$, be an i.i.d sequence with the same distribution as $S_n(I_{n1})$. A similar relation as (6.1) ensures that it suffices to prove that

$$\sum_{i=1}^{k_n} \widetilde{S}_n(I_{ni}) \xrightarrow{d} N(0, \sigma^2(C)). \tag{6.3}$$

We first verify that

$$\operatorname{var}\left(\sum_{i=1}^{k_n} \widetilde{S}_n(I_{ni})\right) = k_n \operatorname{var}(\widetilde{S}_n(I_{n1})) \to \sigma^2(C). \tag{6.4}$$

We have

$$\begin{aligned} k_n \operatorname{var}(\widetilde{S}_n(I_{n1})) &= \operatorname{var}\left(\sum_{i=1}^{m_n} I_t\right) \\ &= m_n \operatorname{var}(I_0) + 2\sum_{h=1}^{m_n-1}(m_n - h)\operatorname{cov}(I_0, I_h). \end{aligned}$$

By regular variation,

$$m_n \operatorname{var}(I_{\{\mathbf{X}_t/a_m \in C\}}) \to \mu(C). \tag{6.5}$$

Fix $k \geq 1$. Then,

$$\left(\sum_{h=1}^{k} + \sum_{h=k+1}^{r_n} + \sum_{h=r_n+1}^{m_n-1}\right)(m_n - h)\operatorname{cov}(I_0, I_h) = R_1 + R_2 + R_3. \tag{6.6}$$

By the same argument as for (6.5), $R_1 \to \sum_{h=1}^{k} \tau_h(C)$ and similar arguments as those for $I_{22}$ and $I_{23}$ in the proof of Theorem 3.1 show that (6.4) holds.

We apply the central limit theorem for the triangular array of i.i.d. mean-zero random variables $\widetilde{S}_n(I_{ni})$, $i = 1, \ldots, k_n$. By Gnedenko and Kolmogorov [17], Theorem 3, page 101, or Theorem 4.1 in Petrov [28], and since (6.4) holds, one needs to verify the following condition for any $\epsilon > 0$:

$$k_n E[(\widetilde{S}_n(I_{n1}))^2 I_{\{|\widetilde{S}_n(I_{n1})| > \epsilon\}}] = E\left(\sum_{t=1}^{m_n} \widetilde{I}_t\right)^2 I_{\{|\widetilde{S}_n(I_{n1})| > \epsilon\}} \to 0. \tag{6.7}$$

A trivial estimate of the right-hand side is given by

$$c m_n^2 P(|\widetilde{S}_n(I_{n1})| > \epsilon) \leq c\epsilon^{-2} m_n^2 \operatorname{var}(\widetilde{S}_n(I_{n1})) = \mathrm{O}(m_n^3/n) = \mathrm{o}(1). \tag{6.8}$$

Next, we show that (6.7) holds under the conditions (3.10). We have, by the Cauchy–Schwarz and Chebyshev inequalities and (6.4), with $J_m = \sum_{t=1}^{m_n} \widetilde{I}_t$, for constants $c > 0$,

$$(E J_m^2 I_{\{|\widetilde{S}_n(I_{n1})| > \epsilon\}})^2 \leq E J_m^4 P(|\widetilde{S}_n(I_{n1})| > \epsilon) \leq c E J_m^4 \operatorname{var}(\widetilde{S}_n(I_{n1})) = E J_m^4 \mathrm{O}(m_n/n). \tag{6.9}$$

We focus on the fourth moment of the partial sum $J_m$, which can be written as

$$E J_m^4 = \sum_{s,t,u,v=1}^{m_n} E\widetilde{I}_{stuv}$$

$$= \sum_{s=1}^{m_n} E\widetilde{I}_s^4 + c_1 \sum_{s\neq t}^{m_n} E(\widetilde{I}_s \widetilde{I}_t^3) + c_2 \sum_{s<t}^{m_n} E(\widetilde{I}_s^2 \widetilde{I}_t^2) + c_3 \sum_{s<t<u<v} E\widetilde{I}_{stuv} + \mathrm{o}(1)$$
$$= A_1 + c_1 A_2 + c_2 A_3 + c_3 A_4 + \mathrm{o}(1).$$

Now,

$$EA_1 = m_n E\widetilde{I}_1^4 \leq m_n E|\widetilde{I}_1| \leq 2m_n p_0 = \mathrm{O}(1).$$

Since $\widetilde{I}_t$ has mean 0, we have

$$\begin{aligned}|A_2| &= \left|\left(\sum_{|s-t|\leq r_n} + \sum_{|s-t|>r_n}\right) E(\widetilde{I}_s \widetilde{I}_t^3)\right| \\ &\leq 2m_n r_n p_0 + \sum_{|s-t|>r_n} |E(\widetilde{I}_s \widetilde{I}_t^3) - E\widetilde{I}_s E\widetilde{I}_t^3)| \\ &\leq 2m_n r_n p_0 + 4m_n \sum_{j=r_n}^{m_n} \alpha_j = \mathrm{O}(r_n).\end{aligned}$$

The third term can be dealt with as follows:

$$\begin{aligned}A_3 &= \left(\sum_{|s-t|\leq r_n} + \sum_{|s-t|>r_n}\right) E(\widetilde{I}_s^2 \widetilde{I}_t^2) \\ &\leq 2m_n r_n p_0 + \sum_{|s-t|>r_n} (|E(\widetilde{I}_s^2 \widetilde{I}_t^2) - E\widetilde{I}_s^2 E\widetilde{I}_t^2| + E\widetilde{I}_s^2 E\widetilde{I}_t^2) \\ &\leq 2m_n r_n p_0 + 4m_n \sum_{j=r_n}^{m_n} \alpha_j + 4m_n^2 p_0^2 = \mathrm{O}(r_n).\end{aligned}$$

We decompose the index set of the fourth term into four disjoint sets:

$$K_1 = \{(s,t,u,v)\colon\ 1 \leq s < t < u < v \leq m_n, v-u > r_n\};$$
$$K_2 = \{(s,t,u,v)\colon\ 1 \leq s < t < u < v \leq m_n, v-u \leq r_n, u-t > r_n\};$$
$$K_3 = \{(s,t,u,v)\colon\ 1 \leq s < t < u < v \leq m_n, v-u \leq r_n, u-t \leq r_n, t-s > r_n\};$$
$$K_4 = \{(s,t,u,v)\colon\ 1 \leq s < t < u < v \leq m_n, v-u \leq r_n, u-t \leq r_n, t-s \leq r_n\}.$$

We then obtain

$$\begin{aligned}|A_4| &= \left|\left(\sum_{K_1} + \sum_{K_2} + \sum_{K_3} + \sum_{K_4}\right) E\widetilde{I}_{stuv}\right| \\ &\leq m_n^3 \sum_{j=r_n}^{m_n} \alpha_j + \left(cm_n^2 r_n \sum_{j=r_n}^{m_n} \alpha_j + \mathrm{O}(r_n^2)\right) + m_n r_n^3 \sum_{j=r_n}^{m_n} \alpha_j + \mathrm{O}(r_n^3)\end{aligned}$$

$$\leq c m_n^3 \sum_{j=r_n}^{m_n} \alpha_j + m_n r_n^3 \sum_{j=r_n}^{m_n} \alpha_j + \mathrm{O}(r_n^3),$$

where the bounds for the sums in the penultimate line follow in the spirit of the arguments used to derive the orders for $A_2$ and $A_3$.

Combining the bounds above with the bound (6.9), the conditions in (3.10) ensure that (6.7) is satisfied. This completes the proof. □

Relation (6.4) suggests the following estimator for $\sigma^2(C)$:

$$\widehat{\sigma}_n^2(C) = k_n^{-1} \sum_{i=1}^{k_n} \left[ \sum_{t \in I_{ni}} I_{\{\mathbf{X}_t/a_m \in C\}} - \frac{m_n}{k_n} \sum_{t=1}^{n} I_{\{\mathbf{X}_t/a_m \in C\}} \right]^2 .$$

It is a consistent estimator of $\sigma^2(C)$, as the following calculations show. We have

$$\widehat{\sigma}_n^2(C) = \sum_{i=1}^{k_n} S_n^2(I_{ni}) - k_n^{-1} S_n^2.$$

Since $k_n^{-1} S_n^2 = \mathrm{o}_P(1)$, it suffices to show that

$$\widetilde{\sigma}_n^2(C) = \sum_{i=1}^{k_n} S_n^2(I_{ni}) \xrightarrow{P} \sigma^2(C).$$

As for (6.2), we observe that for $s \geq 0$,

$$\left| E \prod_{l=1}^{k_n} \mathrm{e}^{-s S_n^2(\widetilde{I}_{nl})} - E \prod_{l=1}^{k_n} \mathrm{e}^{-s \widetilde{S}_n^2(\widetilde{I}_{nl})} \right| \leq 4 k_n \alpha_{r_n} \to 0. \tag{6.10}$$

On the other hand, we have proven that $\sum_{l=1}^{k_n} \widetilde{S}_n(\widetilde{I}_{nl}) \xrightarrow{d} N(0, \sigma^2(C))$. By Raikov's theorem, this is equivalent to

$$\sum_{l=1}^{k_n} \widetilde{S}_n^2(\widetilde{I}_{nl}) \xrightarrow{P} \sigma^2(C).$$

Therefore, from (6.10),

$$\sum_{l=1}^{k_n} S_n^2(\widetilde{I}_{nl}) \xrightarrow{P} \sigma^2(C). \tag{6.11}$$

Let $(W_i)$ be an i.i.d. standard normal sequence, independent of $(\mathbf{X}_t)$. Then, (6.11) holds if and only if

$$W_1^2 \sum_{l=1}^{k_n} S_n^2(\widetilde{I}_{nl}) \stackrel{P}{\to} W_1^2 \sigma^2(C).$$

On the other hand,

$$\begin{aligned} W_1^2 \sum_{l=1}^{k_n} S_n^2(\widetilde{I}_{nl}) &\stackrel{d}{=} \left( \sum_{l=1}^{k_n} W_i S_n(\widetilde{I}_{nl}) \right)^2 \\ &= \left( \sum_{l=1}^{k_n} W_i S_n(I_{nl}) - \sum_{l=1}^{k_n} W_i S_n(J_{ni}) \right)^2 . \end{aligned} \tag{6.12}$$

By virtue of (6.1) ,

$$\sum_{l=1}^{k_n} W_i S_n(J_{ni}) \stackrel{d}{=} W_1 \left( \sum_{i=1}^{k_n} S_n^2(J_{ni}) \right)^{1/2} \stackrel{P}{\to} 0$$

and from (6.12), we conclude

$$\left( \sum_{l=1}^{k_n} W_i S_n(I_{nl}) \right)^2 \stackrel{d}{=} W_1^2 \widetilde{\sigma}_n^2(C) \stackrel{P}{\to} W_1^2 \sigma^2(C),$$

which proves that $\widetilde{\sigma}_n^2(C) \stackrel{P}{\to} \sigma^2(C)$.

## Acknowledgments

We would like to thank Cătălin Stărică for providing the high frequency data set. The comments of both referees were helpful for improving the presentation of the paper. We thank Vladas Pipiras for drawing our attention to the paper by Hill [20]. The first author's research was partly supported by NSF Grants DMS-0338109 and DMS-0743459. The second author's research was partly supported by the Danish Research Council (FNU) Grant 272-06-0442.